\documentclass[12pt, reqno]{amsart}
\usepackage{amssymb,a4}
\usepackage{amsmath,amsthm,amsfonts,amssymb}

\hoffset \voffset \oddsidemargin=60pt \evensidemargin=45pt
\def\fc{{\textbf{\textit c}}}
\def\D{\Delta}
\def\a{\alpha}
\def\b{\beta}

\def\LL{{\mathcal L}}
\def\DD{{\mathcal D}}

\def\Der{{\rm Der}}
\def\Inn{{\rm Inn}}
\def\Ker{{\rm Ker}}
\def\Im{{\rm Im}}

\def\rar{\rightarrow}

\def\sg{\sigma}
\def\vs{\vspace*}

\def\ni{\noindent}
\def\VV{\mathcal {V}}
\def\Z{\mathbb{Z}}
\def\C{\mathbb{C}}
\def\F{\mathbb{C}}
\def\QED{\hfill$\Box$}

\def\ot{\otimes}

\def\pt{\partial}

\def\de{\delta}

\def\a{{\alpha}}

\def\c{{\mathbb C}}

\def\z{{\mathbb Z}}

\def\bN{{\mathbb N}}

\def\b{\beta}

\def\Hom{\mbox{\rm Hom}}

\def\Ker{\mbox{\rm Ker}\,}

\def\Im{\text{\rm Im}}

\numberwithin{equation}{section}
\newtheorem{theo}{Theorem}[section]
\newtheorem{defi}[theo]{Definition}
\newtheorem{coro}[theo]{Corollary}
\newtheorem{lemm}[theo]{Lemma}
\newtheorem{prop}[theo]{Proposition}

\newtheorem{rema}{Remark}

\title{Lie bialgebra structures on the twisted Heisenberg-Virasoro
algebra}
\date{}

\author{Dong Liu}
\address{Department of Mathematics, Huzhou Teachers College, Zhejiang Huzhou, 313000, China}
\email{liudong@hutc.zj.cn}

\author{Yufeng Pei}
\address{Department of Mathematics, Shanghai Normal University, Shanghai, 200234, China}
\email{peiyufeng@gmail.com}

\author{Linsheng Zhu}
\address{Department of Mathematics, Changshu Institute of
Technology, Jiangsu Changshu, 215500, China}
\email{lszhu@cslg.edu.cn}

\begin{document}
\maketitle

\noindent{{\bf Abstract.} In this paper we investigate  Lie
bialgebra structures on the twisted Heisenberg-Virasoro algebra.
With the determination of certain Lie bialgebra structures on the
Virasoro algebra, we determine certain structures on the twisted
Heisenberg -Virasoro algebra. Moreover, some general and useful
results are obtained. With our methods and results we also can
easily to determine certain structures on some Lie algebras related
to the twisted Heisenberg-Virasoro algebra.\vs{5pt}

}

\parskip .001 truein\baselineskip 6pt \lineskip 6pt

\section{Introduction}

\setcounter{section}{1}\setcounter{equation}{0}

Witt type Lie bialgebras were studied in \cite{M, M1, T, NT}, whose
generalized cases were considered in \cite{SS, WSS0}. Lie bialgebra
structures on some Lie (super)algebras including the
Schr\"{o}dinger-Virasoro Lie algebra, the Lie algebra of Weyl type,
the $N=2$ superconformal algebra, were investigated case by case in
\cite{HLS}, \cite{FLL}, \cite{YueS}, \cite{LS}, \cite{YS}, etc..
However, the calculations in these papers are very complicated.
These algebras are all related to the twisted Heisenberg-Virasoro
algebra, which has been first studied by Arbarello et al. in
\cite{ACKP}, where a connection is established between the second
cohomology of certain moduli spaces of curves and the second
cohomology of the Lie algebra of differential operators of order at
most one. Moreover, the twisted Heisenberg-Virasoro algebra has some
relations with the full-toroidal Lie algebras and the $N = 2$ super
conformalalgebra, which is one of the most important algebraic
objects in superstring theory. The structure and representations for
the twisted Heisenberg-Virasoro algebra  was studied in \cite{ACKP},
\cite{B}, \cite{SJ}, \cite{LZ}, \cite{LJ} and \cite{LWZ} etc.
However, Lie bialgebra structures on the twisted Heisenberg-Virasoro
Lie algebra have not yet been considered. Drinfel'd \cite{D2} posed
the problem whether or not there exists a general way to quantize
all Lie bialgebras. Although Etingof and Kazhdan \cite{EK} gave a
positive answer to the question, they did not provide a uniform
method to realize quantizations of all Lie bialgebras. As a matter
of fact, investigating Lie bialgebras and quantizations is a
complicated problem.

In this paper, we shall obtain some Lie bialgebra structures on the
twisted Heisenberg-Virasoro algebra. Moreover we obtain some general
results for this kinds of Lie algebras and provide a uniform method
(no need some complicated calculations) to obtain Lie
(super)bialgebra structures on a series of Lie (super) algebras
related to twisted Heisenberg-Virasoro algebra, including the Lie
algebra of differential operators, the Schr\"{o}dinger-Virasoro
algebra, the $N=2$ superconformal algebra, etc..

We note that the center of the twisted Heisenberg-Virasoro algebra
is 4-dimensional, which is different from the above algebras whose
centers are no more than one-dimensional. So our construction is new
and it has the potential to construct more Lie bialgebra structures.

Throughout the paper, we denote by $\Z_+$ the set of all nonnegative
integers and $\Z^*$ (resp. $\F^*$) the set of all nonzero elements
of $\Z$ (resp. $\F$).

\section{Basics}

\subsection{The twisted Heisenberg-Virasoro algebra}

By definition, as a vector space over $\C$, the twisted
Heisenberg-Virasoro algebra ${\mathcal L}$ has a basis $\{L_m, I_m,
C_L, C_I, C_{LI}\mid m\in\z\}$, subject to the following relations:
\begin{eqnarray}\label{def1}
&&[L_m, L_n]=(n-m)L_{m+n}+\delta_{m+n, 0}{1\over 12}(m^3-m)C_L;\nonumber\\
&&[I_m, I_n]=n\de_{m+n, 0}C_{I};\nonumber\\
&&[L_m, I_n]=nI_{m+n}+\de_{m+n, 0}(m^2-m)C_{LI};\nonumber\\
&&[{\mathcal L}, C_L]=[{\mathcal L}, C_I]=[{\mathcal L}, C_{LI}]=0.
\end{eqnarray}

Clearly the Heisenberg algebra $H=\F\{I_m, C_I\mid m\in\z\}$ and the
Virasoro algebra ${\frak v}=\F\{L_m, C_L\mid m\in\z\}$ are
subalgebras of ${\mathcal L}$. Moreover, ${\mathcal
L}=\oplus_{m\in\z}{\mathcal L}_m$, where ${\mathcal L}_m=\F\{L_m,
I_m\}\oplus\de_{m, 0}\F\{C_I, C_L, C_{LI}\}$, is a $\z$-graded Lie
algebra. Clearly ${\mathcal Z}=\{I_0, C_{LI}, C_L, C_I\}$ is a basis
of the center $\mathcal C$ of $\mathcal L$.

Let ${\mathbb{C}}[t, t^{-1}]$ be the algebra of Laurent polynomials
over $\c$. Denote ${\mathcal G}=\hbox{Diff}\,{\mathbb{C}}[t,
t^{-1}]$ by the Lie algebra of differential operators over
${\mathbb{C}}\,[t, t^{-1}]$. Denote by $D=t\pt$ then as a vector
space over $\c$, ${\mathcal G}=\hbox{Span}_{\c}\{t^mD^n\mid m\in
\mathbb{Z}, \  n\in\bN\}$ with Lie bracket
\begin{equation}
[\,t^mD^n, t^{m_1}D^{n_1}\,]=t^{m+m_1}\left(\sum_{i=1}^n{n\choose
i}m_1^iD^{n+n_1-i}-\sum_{j=1}^{n_1}{n_1\choose
j}m^jD^{n+n_1-j}\right).
\end{equation}

Let ${\mathcal G}_1$ be subalgebra of the Lie algebra ${\mathcal G}$
of differential operators generated by $\{t^m, t^m D\mid m\in\z\}$.
Then ${\mathcal G}_1$ is isomorphic to the centerless twisted
Heisenberg-Virasoro algebra by $t^{m}D$ to $L_m$ and $t^m$ to $I_m$.

Let $M$ be a $ {\mathbb{Z}}$-graded $\LL$-module. Denote by
$\Der(\LL,M)$ the set of \textit{derivations} $\phi:\LL\to M$,
namely, $\phi$ is a linear map satisfying
\begin{eqnarray}
\label{deriv} \phi([x,y])=x\cdot \phi(y)-y\cdot \phi(x),
\end{eqnarray}
and the set $\Inn(\LL,M)$ consisting of the derivations $v_{\rm
inn}, \, v\in M$, where $v_{\rm inn}$ is the \textit{inner
derivation} defined by $v_{\rm inn}:x\mapsto x\cdot v.$ Then it is
well known that $H^1(\LL,M)\cong\Der(\LL,M)/\Inn(\LL,M), $ where
$H^1(\LL,M)$ is the {\it first cohomology group} of the Lie algebra
$\LL$ with coefficients in the $\LL$-module $M$.

 A derivation
$\phi\in\Der(\LL,M)$ is {\it homogeneous of degree $m\in\Z$} if
$\phi(\LL_p)\subset M_{m+p}$ for all $p\in\ \mathbb{Z}$.
 Denote $\Der(\LL,M)_m =\{\phi\in\Der(\LL,M)\,|\,{\rm deg\,}\phi=
m\}$ for $m\in\Z.$ Let $\phi$ be an element of $\Der(\LL,M)$. For
any $m\in\Z$, define the linear map $\phi_m:\LL\rightarrow M$ as
follows: For any $u\in\LL_q$ with $q\in\ \mathbb{Z}$, write
$\phi(u)=\sum_{p\in\ \mathbb{Z}}u_p$ with $u_p\in M_p$, then we set
$\phi_m(u)=u_{q+m}$. Obviously, $\phi_m\in \Der(\LL,M)_m$ and we
have
\begin{eqnarray}\label{summable}
\phi=\mbox{$\sum\limits_{m\in\mathbb{Z}}\phi_m$},
\end{eqnarray}
which holds in the sense that for every $u\in\LL$, only finitely
many $\phi_m(u)\neq 0,$ and $\phi(u)=\sum_{m\in\
\mathbb{Z}}\phi_m(u)$ (we call such a sum in (\ref{summable}) {\it
summable}).

It is well known that ${\mathcal L}$ is the universal central
extension of ${\mathcal G}_1$(see \cite{ACKP}).

\begin{lemm}\cite{SJ}\label{der1}  $H^1({\mathcal L}, {\mathcal L})=\frak D$, where $\frak D$ is consist of the following derivations $\chi$:
\begin{eqnarray*}
&& \chi(L_n)=(\a n+\gamma)I_n+\delta_{n,0}(\gamma+\a)C_{LI},\\
&& \chi(I_n)=\b I_n+\delta_{n,0}(\a+\gamma)C_{I},\\
&& \chi(C_{L})=-24\a C_{LI},\ \chi(C_{LI})=\b C_{LI}-\a C_I,\
\chi(C_{I})=2\b C_I,
\end{eqnarray*} where $\a,\b,\gamma\in \mathbb{C}$, $n\in\z$.
\end{lemm}

The following lemma is very useful to investigate Lie bialgebra
structures for some Lie algebras related to the Virasoro algebra.

\begin{lemm} \label{der}
Suppose that $\frak g=\bigoplus_{n\in\z}{\frak g}_n$ is a
$\z$-graded Lie algebra with a finite dimensional center ${\mathcal
C}_{\frak g}$, and ${\frak g}_0$ is generated by $\{{\frak g}_n,
n\ne 0\}$, then
$$H^1({\frak g},
{\mathcal C}_{\frak g}\ot {\frak g}+{\frak g}\ot {\mathcal C}_{\frak
g})_0={\mathcal C}_{\frak g}\ot H^1({\frak g}, \frak g)_0+H^1({\frak
g}, \frak g)_0\ot {\mathcal C}_{\frak g},$$ where for any $z_1\ot
\sg'+\sg''\ot z_2\in {\mathcal C}_{\frak g}\ot \Der({\frak g},
{\frak g})+\Der({\frak g},{\frak g})\ot {\mathcal C}_{\frak g}$,
$z_1\ot \sg'+\sg''\ot z_2$ is an element of $\Der({\frak g},
{\mathcal C}_{\frak g}\ot {\frak g}+{\frak g}\ot {\mathcal C}_{\frak
g})$ by $(z_1\ot \sg'+\sg''\ot z_2)(x)=z_1\ot \sg'(x)+\sg''(x)\ot
z_2$.
\end{lemm}

\noindent{\bf Proof.} Suppose that $r=\dim {\mathcal C}_{\frak g}$,
and ${\mathcal C}_{\frak g}=\c\{z_1, z_2, \cdots, z_r\}$, then for
any $n\ne 0$, $x\in {\frak g}_n$,
$\sg(x)=\sum_{i=1}^r\sg_i'(x)\otimes z_i+\sum_{i=1}^r z_i\otimes
\sg_i''(x)$ for some $\sg_i', \sg_i''\in\Hom_\c({\frak g}, {\frak
g})$. Then applying $\sg$ to $[x, y]$ for any $x\in {\frak g}_m,
y\in {\frak g}_n$ with $m+n\ne 0$, we have
\begin{eqnarray}
&&\sg_i'([x, y])\ot z_i+z_i\ot \sg_i''([x, y])\nonumber\\
=&& [\sg_i'(x), y]\otimes z_i+[x, \sg_i'(y)]\otimes z_i+z_i\otimes
[\sg_i''(x), y]+z_i\otimes[x, \sg_i''(y)].\nonumber
\end{eqnarray}
Then
\begin{eqnarray*}
&&\sg_i'([x, y])=[\sg_i'(x), y]+[x, \sg_i'(y)],\\
&&\sg_i''([x, y])=[\sg_i''(x), y]+[x, \sg_i''(y)].
\end{eqnarray*}

Since ${\frak g}_0$ is generated by $\{{\frak g}_n, n\ne 0\}$, then
$\sg$ induces derivations $\sg_i', \sg_i''\in\Der({\frak g}, {\frak
g})$ for $i=1, 2\cdots, r$. Moreover $\sg$ is an inner derivation if
and only if all $\sg_i', \sg_i'', i=1, 2, \cdots, r$ are inner
derivations.\qed

\subsection{Lie bialgebras}

Let us recall the definitions related to Lie bialgebras. Let $L$ be
any vector space. Denote $\xi$ the {\it cyclic map} of $L\otimes
L\otimes L$, namely, $ \xi (x_{1} \otimes x_{2} \otimes x_{3})
=x_{2} \otimes x_{3} \otimes x_{1}$ for $x_1,x_2,x_3\in L,$ and
$\tau$ the {\it twist map} of $L\otimes L$, i.e., $\tau(x\otimes y)=
y \otimes x$ for $x,y\in L$. The definitions of a Lie algebra and
Lie coalgebra can be reformulated as follows. A {\it Lie algebra} is
a pair $(L,\delta)$ of a vector space $L$ and a bilinear map $\delta
:L\otimes L\rar L$ with the conditions:
\begin{eqnarray*}
&&\Ker(1-\tau) \subset \Ker\,\delta,\ \ \ \delta \cdot (1 \otimes
\delta ) \cdot (1 + \xi +\xi^{2}) =0 :  L \otimes L\otimes L\rar
 L.
\end{eqnarray*}
Dually, a {\it Lie coalgebra} is a pair $(L,\D)$ of a vector space $
L$ and a linear map $\D: L\to L\otimes L$ satisfying:
\begin{eqnarray}\label{cLie-s-s}
&&\Im\,\D \subset \Im(1- \tau),\ \ \ (1 + \xi +\xi^{2}) \cdot (1
\otimes \D) \cdot \D =0: L\to L\otimes L\otimes L.
\end{eqnarray}
We shall use the symbol ``$\cdot$'' to stand for the {\it diagonal
adjoint action}:
\begin{eqnarray*}
&&x\cdot(\mbox{$\sum\limits_{i}$}{a_{i}\otimes b_{i}})=
\mbox{$\sum\limits_{i}$}({[x,a_{i}]\otimes
b_{i}+a_{i}\otimes[x,b_{i}]}).
\end{eqnarray*}
%\begin{defi}\rm
A {\it Lie bialgebra} is a triple $( L,\delta,\D)$ such that $( L,
\delta)$ is a Lie algebra, $( L,\D)$ is a Lie coalgebra, and the
following compatible condition holds:
\begin{eqnarray}
\label{tr}&&\D\delta (x\otimes y) = x \cdot \D y - y \cdot \D x,\ \
\forall\,\,x,y\in L.
\end{eqnarray}
%\end{defi}
Denote  $U$ the universal enveloping algebra of $ L$, and  $1$ the
identity element of $U$. For any $r =\sum_{i} {a_{i} \otimes
b_{i}}\in L\otimes L$, define $\fc(r)$ to be elements of $U \otimes
U \otimes U$ by
\begin{eqnarray*}
&&\fc(r) = [r^{12} , r^{13}] +[r^{12} , r^{23}] +[r^{13} , r^{23}],
\end{eqnarray*}
where $r^{12}=\sum_{i}{a_{i} \otimes b_{i} \otimes 1} , \ \ r^{13}=
\sum_{i}{a_{i} \otimes 1 \otimes b_{i}} , \ \ r^{23}=\sum_{i}{1
\otimes a_{i} \otimes b_{i}}$. Obviously
\begin{eqnarray*}
\fc(r)=\mbox{$\sum\limits_{i,j}$}[a_i,a_j]\otimes b_i\otimes
b_j+\mbox{$\sum\limits_{i,j}$}a_i\otimes [b_i,a_j]\otimes b_j+
\mbox{$\sum\limits_{i,j}$}a_i\otimes a_j\otimes [b_i,b_j].
\end{eqnarray*}

\begin{defi}\label{def2}
\rm (1) A {\it coboundary Lie bialgebra} is a $4$-tuple $( L,
\delta, \D,r),$ where $( L,\delta,\D)$ is a Lie bialgebra and $r \in
\Im(1-\tau) \subset L\otimes L$ such that $\D=\D_r$ is a {\it
coboundary of $r$}, where $\D_r$ is defined by
\begin{eqnarray}
\label{D-r}\D_r(x)=x\cdot r\mbox{\ \ for\ \ }x\in L.
\end{eqnarray}

(2) A coboundary Lie bialgebra $( L,\delta,\D,r)$ is called {\it
triangular} if it satisfies the following {\it classical Yang-Baxter
Equation} (CYBE):
\begin{eqnarray}
\label{CYBE} \fc(r)=0.
\end{eqnarray}

(3) An element $r\in\Im(1-\tau)\subset L\otimes L$ is said to
satisfy the \textit{modified Yang-Baxter equation} (MCYBE) if
\begin{eqnarray}
\label{MYBE}x\cdot \fc(r)=0,\ \,\forall\,\,x\in L.
\end{eqnarray}
\end{defi}

\begin{lemm}\label{Legr}\rm
Regard $\mathcal L^{\otimes n}$ $($the tensor product of $n$ copies of $\mathcal L)$
as an $\mathcal L$-module under the adjoint diagonal action of $\mathcal L$. Suppose
$r\in \mathcal L^{\otimes n}$ satisfying $x\cdot r=0$, $\forall$ $x\in L$.
Then $r\in {\mathcal C}_{\mathcal L}^{\otimes n}$, where ${\mathcal C}_{\mathcal L}$ is
the center of $\mathcal L$.
\end{lemm}\vs{-6pt}
\noindent{\bf Proof.}\ \ It can be proved directly by using the
similar arguments as those presented in the proof of Lemma 2.2 of
\cite{WSS}(also see the proof Lemma 2.5 of \cite{LCZ}).\QED
\begin{lemm} \label{some}

 Let $L$ be a Lie
algebra and $r\in\Im(1-\tau)\subset L\otimes L,$\  then
\begin{eqnarray}
\label{add-c}(1+\xi+\xi^{2})\cdot(1\otimes\D_r)\cdot\D_r(x)=x\cdot
\fc(r),\ \ \forall\,\,x\in L,
\end{eqnarray}
and the triple $(L,[\cdot,\cdot], \D_r)$ is a Lie bialgebra if and
only if $r$ satisfies MCYBE $(\ref{MYBE})$.
\end{lemm}
\noindent{\bf Proof. }  The result can be found in
\cite{D1,D2,NT}.\QED\vskip5pt

The Lie bialgebra structures over the Witt algebra $W$ (the
centerless Virasoro algebra) and the Virasoro algebra $\frak v$ were
determined in \cite{NT}.

\begin{prop}\cite{NT,SS}\label{P1} For the Witt algebra $W$ and the Virasoro algebra $\frak v$,
$H^1(W, W\otimes W)=H^1({\frak v},{\frak v}\otimes {\frak v})=0$,
and every Lie bialgebra structure on $W$ or $\frak v$ is triangular
coboundary.
\end{prop}

\section{Lie bialgebra structures on the twisted Heisenberg-Virasoro
algebra}

Regard $\VV=\LL\otimes\LL$ as a $\LL$-module under the adjoint
diagonal action,then $\LL$ and $\VV$ are both $
{\mathbb{Z}}$-graded. Now we shall calculate $H^1(\mathcal L, \VV)$
for the twisted Heisenberg-Virasoro algebra $\mathcal L$ with
Proposition \ref{P1}, and then determine Lie bialgebra structures on
the algebra.

For any 6 elements $\a,\a^\dag,\b,\b^\dag,\gamma,\gamma^\dag\in\F$
and $z_1, z_1^\dag, w_1, w_1^\dag\in{\mathcal Z}$, one can easily
verify that the linear map $\varrho:\LL\to\VV$ defined below is a
derivation:

\begin{eqnarray}\varrho(L_n)&=&(n\a+\gamma)z_1\otimes
I_n+(n\a^\dag+\gamma^\dag)I_n\otimes z_1^\dag\nonumber\\
&&+\de_{n, 0}\left((\gamma+\a)z_1\otimes C_{LI}+(\gamma^\dag+\a^\dag)C_{LI}\otimes z_1^\dag\right),\nonumber\\
\varrho(I_n)&=&\b w_1\otimes I_n+\b^\dag I_n\otimes w_1^\dag+\de_{n,
0}\left((\gamma+\a)z_1\otimes
C_{I}+(\gamma^\dag+\a^\dag)C_{I}\otimes z_1^\dag\right),\nonumber\\
\varrho(C_L)&=&-24(\a z_1\otimes C_{LI}+\a^\dag C_{LI}\otimes z_1^\dag),\nonumber\\
\varrho(C_{LI})&=&(\b w_1\otimes C_{LI}+\b^\dag C_{LI}\otimes w_1^\dag)-(\a z_1\otimes C_{I}+\a^\dag C_{I}\otimes z_1^\dag),\nonumber\\
\varrho(C_I)&=&2(\b w_1\otimes C_{I}+\b^\dag C_{I}\otimes
w_1^\dag),\ \ n\in \Z, z_i\in {\mathcal Z}, i=1, 2, 3,
4.\label{def-D}
\end{eqnarray}
Clearly $\varrho$ is an outer derivation of $\Der(\mathcal L
,\mathcal L\otimes\mathcal L)$ if
$\a,\a^\dag,\b,\b^\dag,\gamma,\gamma^\dag$ are not zeros. Denote
$\DD$ the vector space spanned by the such elements $\varrho$ over
$\c$.
 Let $\DD^0$ be the subspace of $\DD$ consisting of
elements $\varrho$ such that $\rho(\LL)\subseteq
\mathrm{Im}(1-\tau)$. Namely, $\DD^0$ is a subspace of $\DD$
consisting of elements $\varrho$ with
$\a=-\a^\dag,\,\b=-\b^\dag,\,\gamma=-\gamma^\dag$, $z_1=z_1^\dag,
w_1=w_1^\dag$.

 The main results of
this paper can be formulated as follows.
\begin{theo}\label{theo}\vskip-3pt
 \begin{itemize}\parskip-2pt
\item[\rm(i)] $\Der(\LL,\VV)=\mathrm{Inn}(\LL,\VV)\oplus\DD$ and
$H^1(\LL,\VV)\cong \DD.$

\item[\rm(ii)]  Let $(\LL,[\cdot,\cdot],\D)$ be a Lie bialgebra such that $\D$
has the decomposition $\D_r+\sg$ with respect to
$\Der(\LL,\VV)=\mathrm{Inn}(\LL,\VV)\oplus\DD$, where
$r\in\VV\,({\rm mod\,}\mathcal C\otimes \mathcal C)$ and $\sg\in
\DD$. Then, $r\in \mathrm{Im}(1-\tau)$ and $\sg\in \DD^0$.
Furthermore, $(\LL,[\cdot,\cdot],\sg)$ is a Lie bialgebra provided
$\sg\in\DD^0$.

\end{itemize}\end{theo}

\noindent{\bf Proof.}  For any $\varphi\in \Der(\LL,\VV)$, we first
claim that if $m\in\Z^*$ then $\varphi_m\in\Inn(\LL,\VV)$. To see
this, denote $\gamma=m^{-1}\varphi_{m}(L_0)\in\VV_{m}$. Then for any
$x_n\in\LL_{n}$, applying $\varphi_{m}$ to $[L_0,x_n]=nx_n$ and
using $\varphi_{m}(x_n)\in \VV_{n+m}$, we obtain
$(m+n)\varphi_{m}(x_n)-x_n\cdot \varphi_{m}(L_0)=L_0\cdot
\varphi_{m}(x_n)-x_n\cdot \varphi_{m}(L_0)=n\varphi_{m}(x_n), $
 i.e.,
$\varphi_{m}(x_n)=\gamma_{\rm inn}(x_n)$. Thus
$\varphi_{m}=\gamma_{\rm inn}$ is inner.

We can claim that $\varphi(L_0)\equiv0( {\rm mod}\, {\mathcal
C}\ot{\mathcal C})$. Indeed, for any $p\in\Z$ and $x_p\in\LL_p$,
applying $\varphi$ to $[L_0,x_p]=px_p$, one has $x_p\cdot
\varphi(L_0)=0( {\rm mod}\, {\mathcal C}\ot{\mathcal C})$. Thus by
Lemma \ref{Legr}, $\varphi(L_0)\equiv0 ( {\rm mod}\, {\mathcal
C}\ot{\mathcal C})$.

Now we claim that for any $\varphi\in{\rm Der}(\LL,\VV)$, (\ref{summable})
is a finite sum. To see this, one can suppose $\varphi_n=(w_n^\dag)_{\rm
inn}$ for some $w_n^\dag\in\VV_n$ and $n\in\Z^*$. If
$\Z'=\{n\in\Z^*\,|\,w_n^\dag\ne0\}$ is an infinite set, then
$\varphi(L_0)=\varphi(L_0)+\sum_{n\in\Z'}L_0\cdot
w_n^\dag=\varphi(L_0)+\sum_{n\in\Z'}nw_n^\dag$ is an infinite sum, which
is not an element in $\VV$, contradicting
 the fact that $\varphi$ is a derivation from $\LL$ to $\VV$. This
 together with Proposition 3.3, 3.4 below
proves  Theorem \ref{theo}(i).

By definition, the algebra ${\frak h}=\C\{I_n, C_I, C_L, C_{LI}\mid
n\in\z\}$ is an ideal of $\mathcal L$. Set $\mathcal H=\mathcal
L\otimes {\frak h}+{\frak h}\otimes \mathcal L$,  then $\mathcal H$
is a $\mathcal L$-submodule of $\VV$. The exact sequence $0\to
\mathcal H\to \mathcal V\to \mathcal V/\mathcal H\to 0$ induces a
long exact sequence
$$\to H^0(\mathcal L, \mathcal K)\to H^1({\mathcal L},{\mathcal H})\to H^1({\mathcal L},
{\mathcal V})\to H^1({\mathcal L},\mathcal K)\to$$ of
$\mathbb{Z}$-graded vector spaces, where all coefficients of the
tensor products are in $\C$, and $\mathcal K={\mathcal V}/\mathcal
H$ is the quotient $\mathcal L$-module, on which $\frak h$ acts
trivially. Clearly $H^0(\mathcal L, \mathcal K)={\mathcal
K}^{\mathcal L}=\{x\in{\mathcal K}\mid {\mathcal L}\cdot x=0\}=0$.
Then $H^1({\mathcal L},{\mathcal H})\cong H^1({\mathcal L},
{\mathcal V})$ if we prove that $H^1({\mathcal L},{\mathcal K})=0$.

\begin{prop}\label{P2}
$H^1({\mathcal L},{\mathcal K})=0.$
\end{prop}

\noindent{\it Proof of Proposition \ref{P2}}.

The exact sequence $0\to \frak h\to {\mathcal L}\to \mathcal L/\frak
h\to 0$ induces an exact sequence
\begin{equation}\label{3eq5}
0\longrightarrow H^1({\mathcal L}/{\frak h},{\mathcal K})
\longrightarrow H^1({\mathcal L},{\mathcal K}) \longrightarrow
H^1({\frak h},{\mathcal K})^{\mathcal L}\end{equation}
 of the 5-term sequence associated to the
Hochschild-Serre spectral sequence  $H^p({\mathcal L}/{\frak h},\\
H^q( {\frak h}, {\mathcal K})\Rightarrow H^{p+q}({\mathcal
L},{\mathcal K})$. Clearly, as $\mathcal L$-modules, the quotient
modules $\mathcal K\cong W\ot W$, on which $\frak h$ acts trivially.
Then $H^1({\mathcal L}/\frak h,{\mathcal K})=0$ by Proposition
\ref{P1}, and $H^1({\frak h},{\mathcal K})^{\mathcal L}$ embeds into
Hom$_{U(W)}(\frak h, W\otimes W)$. So we only need to prove that
Hom$_{U(W)}(\frak h, W\otimes W)=0$.

Let $f \in \text {Hom}_{U(W)}(\frak h, W\otimes W)$, for any $n \in
\mathbb{Z}$,
$$0=f([L_n, C_{LI}])=[L_n, f(C_{LI})].$$ So $f(C_{LI})=0$. Similarly, $f(I_0)=f(C_L)=0$. Moreover, $f(I_m)\in \mathcal V_m$ since
$f([L_0, I_m])=[L_0,f(I_m)]$.

By $f([L_{-m}, I_m])=[L_{-m},f(I_m)]$, we can suppose that $f(I_m)=a_m(L_{2m}\otimes L_{-m}-3L_m\otimes L_0+3L_0\otimes L_m-L_{-m}\otimes L_{2m})$. Moreover by
$mf(I_{n+m})=[L_{n}, f(I_m)]$ for any $n\in\z$ we obtain that $a_{m}=0$ for all
$0\ne m\in\z$. Therefore, $f=0$.  \qed

Now we shall determine $H^1({\mathcal L},{\mathcal H})$. Denote by
${\mathcal L}_C =\mathcal L \otimes \mathcal C+\mathcal C\ot
\mathcal L$, then ${\mathcal L}_C$ is a $\mathcal L$-submodule of
${\mathcal H}$.The exact sequence $0\to {\mathcal L}_C\to \mathcal
H\to \mathcal H/{\mathcal L}_C\to 0$ induces
\begin{equation}\label{eq30}
\to H^0(\mathcal L, \mathcal H/{\mathcal L}_C)\longrightarrow
H^1({\mathcal L},{\mathcal L}_C) \longrightarrow H^1({\mathcal
L},{\mathcal H}) \longrightarrow H^1({\mathcal L},{\mathcal
H}/{\mathcal L}_C)\to\nonumber.
\end{equation}
Clearly $H^0(\mathcal L, \mathcal H/{\mathcal L}_C)=(\mathcal
H/{\mathcal L}_C)^{\mathcal L}=0$. Now we shall prove that
$H^1({\mathcal L},{\mathcal H}/{\mathcal L}_C)=0$, then we have
$H^1({\mathcal L},{\mathcal L}_C)\cong H^1({\mathcal L},{\mathcal
H})$.

\begin{prop}\label{PP3}
$H^1({\mathcal L},\mathcal H/{\mathcal L}_C)=0$.
\end{prop}
\noindent{\it Proof of Proposition \ref{PP3}}.

For any $\varphi\in \Der({\mathcal L},\mathcal H/{\mathcal L}_C)_0$,
$0\ne n\in\Z$, one can write $\varphi(L_n)$ and $\varphi(I_n)$ as
follows
\begin{eqnarray*}
\varphi(L_n)\!&=&\!\mbox{$\sum\limits_{i\in\Z, i\ne
n}$}b_{n,i}L_i\!\otimes\!I_{n-i}\!+\!\mbox{$\sum\limits_{i\in\Z,
i\ne 0}$}b^\dag_{n,i}I_i\!\otimes\! L_{n-i}+\sum_{i\ne0, n}a_{n,
i}I_i\ot I_{n-i},
\end{eqnarray*}
where the sums are all finite, and $a_{n, i}, b_{n,i}, b_{n,
i}^\dag\in\C$ for all $i\in\z$.

For any $n\in\Z$ , the following identities hold,
\begin{eqnarray*}
&&L_1\cdot(L_n\otimes I_{-n})=(n-1)L_{n+1}\otimes
I_{-n}-nL_n\otimes I_{1-n},\\
&&L_1\cdot(I_n\otimes L_{-n})=nI_{n+1}\otimes L_{-n}-(1+n)I_n\otimes
L_{1-n},\\
&&L_1\cdot(I_n\otimes I_{-n})=nI_{n+1}\otimes I_{-n}-nI_n\otimes
I_{1-n}.
\end{eqnarray*}

Let $\triangle$ denote the set consisting of 3 symbols $a,b,b^\dag$.
For each $x\in \triangle$ we define
$I_{x}=\max\{\,|p\,|\,\big|\,x_{1,p}\ne0\}.$ For $n=1$, using the
induction on $\sum_{x\in \triangle} I_x$ in the above identities,
and replacing $\varphi$ by $\varphi-u_{\rm inn}$, where $u$ is a
proper linear combination of $L_{p}\otimes I_{-p}$, $I_{p}\otimes
L_{-p}$ and $I_p\ot I_{-p}$ with $p\in\Z$, one can safely suppose
\begin{eqnarray}\label{B1}
a_{1, k}=b_{1,i}=b^\dag_{1,j}=0\ {\rm for}\ i\neq 0,2,\,j\neq\pm1,
k\in\z.
\end{eqnarray}
Applying $\varphi$ to $[\,L_{1},L_{-1}]=-2L_0$ and using the fact
that $\varphi(L_0)\in {\mathcal C}\ot{\mathcal C}$,  and comparing
the coefficients of $L_p\otimes I_{-p}$, $I_p\otimes L_{-p}$ and
$I_p\ot I_{-p}$, we obtain \vspace{-5pt}
\begin{eqnarray*}
&&\mbox{$\sum\limits_{p\in
\Z}$}\big((p-2)b_{-1,p-1}-(1+p)b_{-1,p}+(p-1)b_{1,p}-(p+2)b_{1,1+p}\big)=0,\\
&&\mbox{$\sum\limits_{p\in
\Z}$}\big((p-1)b^\dag_{-1,p-1}-(p+2)b^\dag_{-1,p}
+(p-2)b^\dag_{1,p}-(p+1)b^\dag_{1,1+p}\big)=0,\\
&&\mbox{$\sum\limits_{p\in \Z}$}(p-1)c_{-1, p-1}-(p-1)c_{-1,
p}+(p-1)c_{1,p}-(p+1)c_{1,p+1}=0.
\end{eqnarray*}
Then
\begin{equation}\label{B2}
b_{\pm1,q}=b_{\pm1,q}^\dag=a_{-1, q}=0,\quad \forall q\in\Z.
\end{equation}
Applying $\varphi$ to $[\,L_{2},L_{-1}]=-3L_{1}$, and comparing the
coefficient of $L_p\otimes I_{1-p}$, $I_p\otimes L_{1-p}$ and
$I_p\ot I_{1-p}$ one can obtain that
\begin{equation}\label{B3}
a_{2, p}=b_{2,p}=b^\dag_{2,p}=0.
\end{equation}

Similarly by $[L_1, L_{-2}]=-3L_{-1}$, we have $a_{-2,p}=0$ for all
$p\in\z$. It follows from this formula and (\ref{B1})--(\ref{B3})
that $\varphi(L_{\pm1})=\varphi(L_{\pm2})=0$. Thus for any $n\in
\Z$, one can deduce $\varphi(L_n)=0$, since $\frak v$ can be
generated by $L_{\pm1}$ and $L_{\pm2}$. While for $I_n$, by
$[L_{-n}, I_{n}]=nI_0$ and $[L_m, I_n]=nI_{m+n}$ we have
$\varphi(I_n)=0$.

\begin{rema}
The Proposition \ref{P2}, \ref{PP3} hold for many Lie algebras
related to the Virasoro algebra. For example, $W(a, b)=\C\{L_m,
I_n\mid m, n\in\z\}$, where $L_m, I_m, m\in\z$ are the centerless
Virasoro, Heisenberg operators, and the twisted action is given by
$[L_m, I_n]=(a+bm +n)I_{m+n}$ for all $m,n\in\Z$ and for some
$a,b\in\C$ see \cite{GJP} for detail.
\end{rema}

\begin{prop}\label{P3}
$H^1({\mathcal L},{\mathcal L}_C)=\mathcal D$, where ${\mathcal
L}_C={\mathcal C}\ot {\mathcal L}+{\mathcal L}\ot{\mathcal C}$.
\end{prop}
\noindent{\it Proof of Proposition \ref{P3}}.

For any $\varphi\in \Der({\mathcal L}, {\mathcal L}_C)_0$,  $n\ne
0$, from Lemma \ref{der1}, \ref{der}, we suppose that

\begin{eqnarray}
\varphi(L_n)\!&=&(n\a+\gamma)z_1\otimes I_n\! +\!(n\a^\dag+\gamma^\dag)I_n\!\otimes\!z_1^\dag,\nonumber\\
\varphi(I_n)\!&=&\b w_1\otimes I_n \! +\!\b^\dag
I_n\!\otimes\!w_1^\dag,\label{3eq2}\end{eqnarray} for some $\a,
\a^\dag, \b, \b^\dag, \gamma, \gamma^\dag\in\c$.

By $[L_1, L_{-1}]=-2L_0$ we have
\begin{equation}
\varphi(L_0)=\gamma z_1\otimes I_0+\gamma^\dag I_0\otimes
z_1^\dag+(\a+\gamma)z_1\otimes
C_{LI}+(\a^\dag+\gamma^\dag)C_{LI}\otimes z_1^\dag.\label{3eq3}
\end{equation}

By $[L_{1}, I_{-1}]=-I_0$ we have
\begin{equation}
\varphi(I_0)=\b w_1\otimes I_0+\b^\dag I_0\otimes
w_1^\dag+(\a+\gamma)z_1\otimes
C_{I}+(\a^\dag+\gamma^\dag)C_{I}\otimes z_1^\dag.\label{3eq3}
\end{equation}

Moreover $C_L$, $\varphi(C_{LI})$ and $\varphi(C_I)$ can be computed
by $[L_2, L_{-2}]=-4L_0+{1\over2}C_L$, $[L_{-1}, I_{1}]=I_0+2C_{LI}$
and $[I_1, I_{-1}]=-C_I$. Then we get $\varphi\in\mathcal D$. \qed

Then we get Theorem \ref{theo} (i).

To prove the second part of Theorem 3.1 (ii), we need the following
lemma.
\begin{lemm}\rm \label{lemma3ll}
Suppose $v\in\VV$ such that $x\cdot v\in {\rm Im}(1-\tau)$ for all
$x\in\LL.$ Then there exists $u\in {\rm Im}(1-\tau)$ such that
$v-u\in \mathcal C\otimes \mathcal C$.
\end{lemm}
\ni{\bf Proof}\ \ First note that $\LL\cdot {\rm Im}(1-\tau)\subset
{\rm Im}(1-\tau).$  We prove that after several steps, by replacing
$v$ with $v-u$ for some $u\in {\rm Im}(1-\tau)$, we get $v\in
\mathcal C\otimes \mathcal C$. Write $v=\sum_{n\in\Z}v_n,
v_n\in{\mathcal V}_n$. Obviously,
\begin{eqnarray}\label{eqrx}
v\in {\rm Im}(1-\tau)\ \,\Longleftrightarrow \ \,v_n\in {\rm
Im}(1-\tau),\ \ \forall\,\,n\in\Z.
\end{eqnarray}
Then $\sum_{n\in\Z}nv_n=L_0\cdot v\in {\rm Im}(1-\tau)$. By
(\ref{eqrx}), $nv_n\in {\rm Im}(1-\tau),$ in particular, $v_{n}\in
{\rm Im}(1-\tau)$ if $n\ne0$. Thus by replacing $v$ by
$v-\sum_{n\in\Z^*}v_n$, one can suppose $v=v_0\in\VV_0$. Write
\begin{eqnarray*}
v\!\!\!&=&\!\!\!\mbox{$\sum\limits_{i\in\Z}$}a_{i}L_i\otimes L_{-i}+
\mbox{$\sum\limits_{0\ne p\in\Z}$}(b_{p}L_p\otimes
I_{-p}+c_{p}I_p\otimes L_{-p}+d_{p}I_p\otimes I_{-p})\\
&+&bL_0\otimes z+b^\dag z^\dag\otimes L_0(\hskip-10pt\mod \mathcal
C\otimes \mathcal C),
\end{eqnarray*}
where all the coefficients are in $\F$, $z, z^\dag\in\mathcal Z$ and
the sums are all finite. Since the elements of the form
$u_{1,p}:=L_p\otimes L_{-p}-L_{-p}\otimes L_{p}$,
$u_{2,p}:=L_p\otimes I_{-p}-I_{-p}\otimes L_{p},$
$u_{3,p}:=I_p\otimes I_{-p}-I_{-p}\otimes I_{p}$, $u=L_0\otimes
z-z\otimes L_0$ are all in ${\rm Im}(1-\tau),$ replacing $v$ by
$v-u$, where $u$ is a combination of some $u_{1,p}$, $u_{2,p}$ and
$u_{3,p}$, one can suppose
\begin{eqnarray}\label{wpqr2}
&&c_p=0,\ \forall\ \,p\in\Z;\ \ a_{p},\ d_{p}\ne 0\,\Longrightarrow\
\,p>0\ \mbox{ or }\ p=0.
\end{eqnarray}
Then $v$ can be rewritten as
\begin{eqnarray}\label{sm1}
v&=&\mbox{$\sum\limits_{p\in\Z_{+}}a_{p}$}L_p\otimes
L_{-p}+\mbox{$\sum\limits_{0\ne p\in\Z}b_{p}$}L_p\otimes
I_{-p}+\mbox{$\sum\limits_{p\in\Z_{+}}d_{p}$}I_p\otimes
I_{-p}\nonumber\\
&+&bL_0\otimes z(\hskip-10pt\mod \mathcal C\otimes \mathcal C).
\end{eqnarray}
Assume $a_{p}\ne 0$ for some $p>0$. Choose $q>0$ such that $q\ne p$.
Then $L_{p+q}\otimes L_{-p}$ appears in $L_{q}\cdot v,$ but
(\ref{wpqr2}) implies that the term $L_{-p}\otimes L_{p+q}$ does not
appear in $L_q\cdot v$, a contradiction with the fact that $L_q\cdot
v\in {\rm Im}(1-\tau)$. Then one can suppose $a_p=0,\
\forall\,\,p\in\Z^*$ . Similarly, one can also suppose $d_p=0,\
\forall\,\,p\in\Z^*$ and $e_p=0$, $\forall\,\,p\in\Z$. Then
(\ref{sm1}) becomes
\begin{eqnarray}\label{sm2}
v=\mbox{$\sum\limits_{0\ne p\in\Z}$}b_{p}L_p\otimes I_{-p}
+a_{0}L_0\otimes L_{0}+bL_0\otimes z\mod \mathcal C\otimes \mathcal C.
\end{eqnarray}
Recall the fact ${\rm Im}(1-\tau)\subset{\rm Ker}(1+\tau)$ and our
hypothesis $\LL\cdot v\subset{\rm Im}(1-\tau)$, one has
\begin{eqnarray*}
0\!\!\!&=&\!\!\!(1+\tau)L_1\cdot v\\
\!\!\!&=&\!\!\!-2a_{0}(L_1\otimes L_{0}+L_0\otimes
L_{1})-b(L_1\otimes z+z\otimes L_1)\\
&&\!\!\!+\mbox{$\sum\limits_{0\ne p\in\Z}$}\big((p-1)b_{p}
L_{p+1}\otimes I_{-p}-pb_pL_p\otimes
I_{1-p}\big)+\big((p-1)b_{p}I_{-p}\otimes L_{p+1}
-pb_{p}I_{1-p}\otimes L_p\big).
\end{eqnarray*}
Comparing the coefficients, and noting that the set
$\{p\,|\,b_p\ne0\}$ is finite, one gets
\begin{eqnarray*}
&&a_0=b_{p}=0, p\ne 1.
\end{eqnarray*}

Moreover $b_1=b=0$ if $z\ne kI_0$ for some $k\in\c$, $b+b_1=0$ if
$z=I_0$.

Then (\ref{sm2}) can be rewritten as
\begin{eqnarray}\label{sm3}
v=b_{1}(L_{1}\otimes I_{-1}-L_0\otimes I_{0})(\hskip-10pt\mod
\mathcal C\otimes \mathcal C).
\end{eqnarray}
Observing $(1+\tau)L_2\cdot v=0$, one has $b_1=0$.  Thus the lemma
follows.

\ni{\it Proof of Theorem \ref{theo} (ii)}\ \ \rm\  Let $(\LL
,[\cdot,\cdot],\D)$ be a Lie bialgebra structure on $\LL$.  By
(\ref{deriv}), (\ref{tr}) and Theorem \ref{theo}(i), $\D=\D_r+\sg$,
where $r\in \VV\,({\rm mod\,}\mathcal C\otimes \mathcal C)$ and
$\sg\in\DD$. By (\ref{cLie-s-s}), ${\rm Im}\,\D\subset{\rm
Im}(1-\tau)$, so $\D_r(L_n)+\sg(L_n)\in$ Im$(1-\tau)$ for $n\in\Z$,
which implies that $\a+\a^\dag=\gamma+\gamma^\dag=0$. Similarly,
$\b+\b^\dag=0$ by the fact that
$\D_r(I_n)+\sg(I_n)\in\mathrm{Im}(1-\tau)$ for $n\in\Z$. Thus,
$\sg(\LL)\in \mathrm{Im}(1-\tau)$. So
$\mathrm{Im\,}\D_r\in\mathrm{Im}(1-\tau)$. It follows immediately
from  Lemma \ref{lemma3ll} that $r\in{\rm Im}(1-\tau)\,({\rm
mod\,}\mathcal C\otimes \mathcal C)$, proving the first statement of
Theorem \ref{theo}(ii). If $\sg\in\DD^0$, one can easily verify that
$(1+\xi+\xi^2)\cdot(1\otimes \sg)\cdot \sg=0$ by acting it on
generators of $\LL$, which shows $(\LL,[\cdot,\cdot],\sg)$ is a Lie
bialgebra, and the proof of Theorem \ref{theo}(ii) is completed.
 \QED

\section{Lie bialgebra structures on some other Lie algebras}

With the methods and results in Section 3, we can easily to obtain
Lie bialgebra structures on some other Lie algebras related the
twisted Heisenberg-Virasoro algebra. Although some of them were
studied case by case in some papers (see \cite{LS},\cite{HLS},
\cite{FLL},\cite{LCZ}, etc.), their calculations are very
complicated.

\subsection{Lie algebra of differential operators of order at most one} Let ${\mathcal
G}_1$ be the Lie algebra of differential operators of order at most
one. Then ${\mathcal G}_1$ is just the centerless twisted
Heisenberg-Virasoro algebra.

As a vector space over $\C$, ${\mathcal G}_1$ has a basis $\{L_m,
I_m, m\in\z\}$, subject to the following relations:
\begin{eqnarray*}
&&[L_m, L_n]=(n-m)L_{m+n},\\
&&[I_m, I_n]=0,\\
&&[L_m, I_n]=nI_{m+n}.\\
\end{eqnarray*}

\begin{lemm}\rm \label{some1} For the Lie algebra ${\mathcal G}_1$ and $r\in {\mathcal G}_1\ot {\mathcal G}_1$,
$r$ satisfies CYBE in $(\ref{CYBE})$ if and only if it satisfies
MCYBE in $(\ref{MYBE})$.
\end{lemm}

\begin{rema}  Let $L$ be a Lie algebra such that $H^0(L; \bigwedge^3L) =0$.
Then any solution $r\in L\bigwedge L$ of the MCYBE is actually a
solution of CYBE. All the following Lie algebras satisfy this
condition, then Lemma \ref{some1} holds for all the following Lie
algebras.

\end{rema}

For any 6 elements $\a,\a^\dag,\b,\b^\dag,\gamma,\gamma^\dag\in\F$,
one can easily verify that the linear map $\varrho:\LL\to\VV$
defined below is a derivation:
\begin{eqnarray}\label{def-D}&&\varrho(L_n)=(n\a+\gamma)I_0\otimes
I_n+(n\a^\dag+\gamma^\dag)I_n\otimes I_0,\nonumber\\
&&\varrho(I_n)=\b I_0\otimes I_n+\b^\dag I_n\otimes I_0.
\end{eqnarray}
 Denote $\DD_1$ the vector space spanned by the such elements $\varrho$ over $\c$.
 Let $\DD_1^0$ be the subspace of $\DD_1$ consisting of
elements $\varrho$ such that $D(\LL)\subseteq \mathrm{Im}(1-\tau)$.
Namely, $\DD_1^0$ is a subspace of $\DD_1$ consisting of elements
$\varrho$ with $\a=-\a^\dag,\,\b=-\b^\dag,\,\gamma=-\gamma^\dag$.

 According to the results in Section 3, we obtain the following
 result.
\begin{theo}\label{theo2}

\begin{itemize}\parskip-2pt

\item[\rm(i)] $\Der({\mathcal G}_1,{\mathcal G}_1\ot {\mathcal G}_1)=\mathrm{Inn}({\mathcal G}_1,{\mathcal G}_1\ot {\mathcal G}_1)\oplus\DD_1$ and
$H^1({\mathcal G}_1,{\mathcal G}_1\ot {\mathcal G}_1)=\Der({\mathcal G}_1,{\mathcal G}_1\ot {\mathcal G}_1)/\Inn({\mathcal G}_1,{\mathcal G}_1\ot {\mathcal G}_1)\cong \DD_1.$

\item[\rm(ii)]  Let $({\mathcal G}_1,[\cdot,\cdot],\D)$ be a Lie bialgebra such that $\D$
has the decomposition $\D_r+\sg$ with respect to
$\Der({\mathcal G}_1,{\mathcal G}_1\ot {\mathcal G}_1)=\mathrm{Inn}({\mathcal G}_1,{\mathcal G}_1\ot {\mathcal G}_1)\oplus\DD_1$, where
$r\in{\mathcal G}_1\ot {\mathcal G}_1\,({\rm mod\,}\mathcal C\otimes \mathcal C)$ and $\sg\in \DD_1$. Then,
$r\in \mathrm{Im}(1-\tau)$ and $\sg\in \DD_1^0$. Furthermore,
$(\LL,[\cdot,\cdot],\sg)$ is a Lie bialgebra provided $\sg\in\DD_1^0$.

\item[\rm(iii)] A Lie bialgebra $({\mathcal G}_1,[\cdot,\cdot],\D)$  is triangular coboundary
if and only if $\D$ is an inner derivation $($thus $\D=\D_r$, where
$r\in\mathrm{Im}(1-\tau)$ is a solution of CYBE$)$.
\end{itemize}\end{theo}

\noindent{\bf Proof.}
\ref{theo2}(iii) follows immediately from (\ref{cLie-s-s}),
Definition \ref{def2} and Lemma \ref{some1}.

\subsection{The Lie algebra of differential operators}

Now we consider Lie bialgebra structures on the Lie algebra of
differential operators defined in Section 2. Although such works for the Lie algebra of Weyl type (including the centerless Lie algebra of differential operators) was consider in \cite{YueS}, our calculation is very simple. Clearly ${\mathcal
C}_d=\c\{1\}$ is the center of $\mathcal G$. Moreover, we have

\begin{lemm}\label{lemm0}
 $\mathcal G$ is generated by $\{t, t^{-1}, D^2\}$(see \cite{Z}).
\end{lemm}

For any 2 elements $\a,\a^\dag\in\F$, one can easily verify that the
linear map $\zeta_{\a, \a^\dag}:{\mathcal G}\to {\mathcal G}\ot
{\mathcal G}$ defined below is a derivation\vs{-7pt}:
\begin{eqnarray}&&\zeta_{\a,\a^\dag}(t^m)=0,\nonumber\\
&&\zeta_{\a,\a^\dag}(t^mD^n)=\a nt^mD^{n-1}\otimes 1+\a^\dag n\cdot
1\otimes t^mD^{n-1}, \  \forall m\in\z,
n\in\z_+.\label{def-D}\end{eqnarray}

Clearly $\zeta_{\a, \a^\dag}$ is an outer derivation of
$\Der({\mathcal G} ,{\mathcal G}\ot {\mathcal G})$ if $\a,\a^\dag$
are not zeros. Denote $\DD_d$ the vector space spanned by the such
elements $\zeta_{\a,\a^\dag}$ over $\c$.
 Let $\DD^0_d$ be the subspace of $\DD_d$ consisting of
elements $\zeta_{\a,\a^\dag}$ such that
$\zeta_{\a,\a^\dag}({\mathcal G})\subseteq \mathrm{Im}(1-\tau)$.
Namely, $\DD^0_d$ is a subspace of $\DD_d$ consisting of elements
$\sg$ with $\a=-\a^\dag$.

With the results in Section 3, we can obtain the main results of this
paper as follows.
\begin{theo}\label{theo4}

\rm(i) $\Der({\mathcal G},{\mathcal G}\ot {\mathcal G})
=\mathrm{Inn}({\mathcal G},{\mathcal G}\ot {\mathcal G})\oplus\DD_d$
and $H^1({\mathcal G},{\mathcal G}\ot {\mathcal G})\cong \DD_d.$

\rm(ii)  Let $({\mathcal G},[\cdot,\cdot],\sg)$ be a Lie bialgebra
such that $\sg$ has the decomposition $\D_r+\sg$ with respect to
$\Der({\mathcal G},{\mathcal G}\ot {\mathcal G})=
\mathrm{Inn}({\mathcal G},{\mathcal G}\ot {\mathcal G})\oplus\DD_d$,
where $r\in{\mathcal G}\ot {\mathcal G}\,({\rm mod\,}{\mathcal
C}_d\otimes {\mathcal C}_d)$ and $\sg\in \DD_d$. Then, $r\in
\mathrm{Im}(1-\tau)$ and $\sg\in \DD^0_d$. Furthermore, $({\mathcal
G},[\cdot,\cdot],\sg)$ is a Lie bialgebra provided $\sg\in\DD^0_d$.

\rm(iii) A Lie bialgebra $({\mathcal G},[\cdot,\cdot],\D)$  is
triangular coboundary if and only if $\D$ is an inner derivation
$($thus $\D=\D_r$, where $r\in\mathrm{Im}(1-\tau)$ is some solutions
of CYBE$)$.
\end{theo}
\noindent{\bf Proof.} For any $\sg\in\Der({\mathcal G},{\mathcal
G}\ot {\mathcal G})$, we can suppose that $\sg(t)=0$(see
\cite{YueS}, or \cite{Z}). With $[t^m,t]=0$ and $[t^mD, t]=t^{m+1}$,
we have $\sg(t^m), \sg(t^mD)\in {\mathcal G}_1\ot{\mathcal G}_1$.
Then we have
\begin{equation}
\Der({\mathcal G}_1,{\mathcal G}\ot {\mathcal G})\subset
\Der({\mathcal G}_1,{\mathcal G}_1\ot {\mathcal G}_1)+\Inn({\mathcal
G}_1,{\mathcal G}\ot {\mathcal G})
\end{equation}
By Theorem \ref{theo2}, and replaced $\sg$ by
$\sg-\a\zeta_{1,0}-\a^\dag\zeta_{0,1}$, one can suppose that
\begin{eqnarray}\label{def-D4}&&\sg(t^mD)=\gamma1\otimes
t^m+\gamma^\dag t^m\otimes 1,\nonumber\\
&&\sg(t^m)=0.
\end{eqnarray}

By $[D^2, t^m]=m^2t^m+2mt^mD, [D^2, D]=0$ we have $[\sg(D^2), t^m]=2m\sg(t^mD)$ for all $m\in\z$, and
$[\sg(D^2), D]=0$. Then can suppose that
\begin{equation}\label{3eq21}
\sg(D^2)=2\gamma 1\ot
D+2\gamma^\dag D\ot 1+\sum f_i t^i\ot t^{-i}.\end{equation}

From $[D^2, t^mD]=m^2t^mD+2mt^mD^2$, we have
\begin{eqnarray}
\sg(t^mD^2)&=&2\gamma 1\ot t^mD+2\gamma^\dag t^mD\ot1 \nonumber\\
&+&{1\over 2m}(\sum_{i\ne0}f_i(it^i\ot t^{m-i}-it^{m+i}\ot
t^i)),\end{eqnarray} for $m\ne 0$.

Applying $\sg$ to $[t^{n}D^2,t^{-n}D]=-3nD^2+n^2D$, and comparing
the coefficients of the term $t^i\otimes t^{n-i}$ with lowest first
degree, we get $f_i=0$ if $i\ne0$. Moreover,
$-3f_0n+n^2(\gamma+\gamma^\dag)=0$ holds for all $n\ne0$. Then
$f_0=0$ and $\gamma+\gamma^\dag=0$. Therefore
\begin{eqnarray}
\sg(t^mD^2)=2\gamma t^mD\ot 1+2\gamma^\dag 1\ot t^mD\end{eqnarray} for all $
m\in\z$.

Applying $\sg$ to $[t^{m}D^2,t^nD]=(2n-m)t^{m+n}D^2+n^2t^{m+n}D$,
and comparing the coefficients of terms $t^{m+n}\ot 1$ and $1\ot
t^{m+n}$, we get $\gamma=\gamma^\dag=0$. Therefore $\sg=0$ by Lemma
\ref{lemm0}. Then we get Theorem \ref{theo4} (i).
Theorem \ref{theo4} (ii) and (iii) are easily proved.

\subsection{The Schr\"{o}dinger-Virasoro Lie algebra}

The Schr\"{o}dinger-Virasoro algebra $\mathfrak{sv},$ playing
important roles in mathematics and statistical physics, is a
infinite-dimensional Lie algebra first introduced by M. Henkle in
[6] by looking at the invariance of the free Schr\"{o}dinger
equation in (1+1)dimensions: $(2\mathcal
{M}\partial_t-\partial_r^{2})\psi=0.$ The structure and the
representation theory for it have been well studied by many authors
(see \cite{HLS}, \cite{LS1}, etc.).

The {\it Schr\"{o}dinger-Virasoro Lie algebra} $\mathfrak{sv}$ is
the infinite-dimensional Lie algebra $\mathfrak{sv}$ with
$\C$-basis $\{L_n,I_n, Y_r\,|n\in \Z,r\in\Z+{1\over2}\}$ and the
following relations
\begin{eqnarray}
\!\!\!&\!\!\!&[L_m, L_n]=(n-m)L_{m+n}+\delta_{m+n, 0}{1\over
12}(m^3-m)C,\nonumber
\\[4pt]\!\!\!&\!\!\!&
[L_n,Y_r]=(r-\frac{n}{2})Y_{n+r},\nonumber\\[4pt]\!\!\!&\!\!\!&
[L_m,I_n]=nI_{m+n},%\label{eq1.1}
\nonumber\\[4pt] \!\!\!&\!\!\!&
[Y_r,Y_{s}]=(s-r)I_{r+s},
\nonumber\\[4pt]\!\!\!&\!\!\!&
[Y_r,I_p]=[I_n,I_p]=0.\nonumber
\end{eqnarray}
where $m,n\in\Z,r,s\in\Z+{1\over2}$.

Clearly ${\mathcal C}_s=\C\{I_0,  C\}$ is the center of
$\mathfrak{sv}$. Denote ${\mathcal Z}_s$ by the set $\{I_0,  C\}$.

Lie bialgebra structures over the cenreless Schr\"{o}dinger-Virasoro
Lie algebra ($C=0$) were determined in \cite{HLS}. However, using
the methods and results in Section 3, the original proof in
\cite{HLS} can be greatly simplified.

For any 6 elements $\a,a^\dag,\b,\b^\dag,\gamma,\gamma^\dag\in\F$,
$z_1, z_1^\dag, w_1, w_1^\dag\in{\mathcal Z}_s$, one can easily
verify that the linear map
$\rho:\mathfrak{sv}\to\mathfrak{sv}\otimes\mathfrak{sv}$ defined
below is a derivation:

\begin{eqnarray}\label{od4}
& &\rho(L_n)=(n\a+\gamma)z_1\otimes
I_n+(n\a^\dag+\gamma^\dag))I_n\otimes z_1^\dag,\nonumber\\
& & \rho(I_n)=2(\b w_1\otimes I_n+\b^\dag I_n\otimes w_1^\dag),\nonumber\\
& &\rho(Y_{n-\frac12})=\b w_1\otimes Y_{n-\frac12}+\b^ \dag
Y_{n-\frac12}\otimes w_1^\dag,\nonumber\\
& &\rho(C)=0,\ \ n\in \Z.
\end{eqnarray}

Denote ${\mathcal D}_s$ the vector space spanned by the such
elements $\rho$. Let ${\mathcal D}_s^0$ be the subspace of
${\mathcal D}_s$ consisting of elements $\rho$ such that
$\rho(\mathfrak{sv})\subseteq \mathrm{Im}(1-\tau)$. Namely,
${\mathcal D}_s^0$ is the subspace of ${\mathcal D}_s$ consisting of
elements $\rho$ with
$\a=-\a^\dag,\,\b=-\b^\dag,\,\gamma=-\gamma^\dag$ and $z_1=z_1^\dag,
\, w_1=w_1^\dag$.

\begin{theo}\label{theo3}$($The centerless case was given in
\cite{HLS}$)$

\rm(i)
$\Der(\mathfrak{sv},\mathfrak{sv}\otimes\mathfrak{sv})=\mathrm{Inn}(\mathfrak{sv},\mathfrak{sv}\otimes\mathfrak{sv})\oplus{\mathcal
D}_s$ and
$H^1(\mathfrak{sv},\mathfrak{sv}\otimes\mathfrak{sv})=\Der(\mathfrak{sv},\mathfrak{sv}\otimes\mathfrak{sv})/\Inn(\mathfrak{sv},\mathfrak{sv}\otimes\mathfrak{sv})\cong
{\mathcal D}_s.$

\rm(ii)  Let $(\mathfrak{sv},[\cdot,\cdot],\D)$ be a Lie bialgebra
such that $\D$ has the decomposition $\D_r+\sg$ with respect to
$\Der(\mathfrak{sv},\mathfrak{sv}\otimes\mathfrak{sv})=\mathrm{Inn}(\mathfrak{sv},\mathfrak{sv}\otimes\mathfrak{sv})\oplus{\mathcal
D}_s$, where $r\in\mathfrak{sv}\otimes\mathfrak{sv}\,({\rm
mod\,}\mathcal C\otimes \mathcal C)$ and $\sg\in {\mathcal D}_s$.
Then, $r\in \mathrm{Im}(1-\tau)$ and $D\in {\mathcal D}_s^0$.
Furthermore, $(\mathfrak{sv},[\cdot,\cdot],\sg)$ is a Lie bialgebra
provided $\sg\in{\mathcal D}_s^0$.

\rm(iii)  A Lie bialgebra $(\mathfrak{sv},[\cdot,\cdot],\D)$  is
triangular coboundary if and only if $\D$ is an inner derivation
$($thus $\D=\D_r$, where $r\in\mathrm{Im}(1-\tau)$ is a solution of
CYBE$)$.
\end{theo}

\noindent{\bf Proof.}
Set ${\mathcal L}$ be the subalgebra of $\mathfrak{sv}$ generated
by $\{L_m, I_m, C\}$, then ${\mathcal L}$ is the twisted
Heisenberg-Virasoro Lie algebra with $C_{LI}=C_I=0, C_L=C$. From the
proof of Theorem \ref{theo}(i), we see that the most difficulty of
the proof (i) is to determine the terms of $\varphi(L_n),
\varphi(I_n)$ for $\varphi\in\Der(\mathfrak{sv},
\mathfrak{sv}\otimes \mathfrak{sv})_0$. It is easy to prove that
$\varphi(L_0)\in {\mathcal C}_s\ot{\mathcal C}_s$ (see the proof of
Theorem \ref{theo}(i)).

Now we suppose that
\begin{eqnarray*}
\varphi(L_n)\!&=&\!\mbox{$\sum\limits_{i\in\Z}$}a_{n,i}Y_{i-\frac12}\!\otimes\!
Y_{n-i+\frac12}(\hskip-10pt\mod {\mathcal L}\ot {\mathcal L}),\\
\varphi(I_n)\!&=&\!\mbox{$\sum\limits_{i\in\Z}$}b_{n,i}Y_{i-\frac12}\!\otimes\!
Y_{n-i+\frac12}(\hskip-10pt\mod {\mathcal L}\ot {\mathcal L}),
\end{eqnarray*}
where the sums are all finite, and $a_{n, i},b_{n,i}\in\C$ for all
$i\in\z$.

\begin{rema}
Although ${\mathcal L}$ is not an ideal of $\mathfrak{sv}$, the
terms $Y_{r}\!\otimes\! Y_{s}$ are only obtained from $Y\ot Y$,
where $Y=\C\{Y_r\mid r\in\z+\frac12\}$ by adjoint actions. So we can use the notation
$x=y\; (\hskip-5pt\mod {\mathcal L}\ot {\mathcal L})$ if $x-y\in
{\mathcal L}\ot {\mathcal L}$.
\end{rema}

With \begin{eqnarray*} L_1\cdot(Y_{n-\frac12}\otimes
Y_{\frac12-n})=(n-1)Y_{n+\frac12}\otimes
 Y_{\frac12-n}-nY_{n-\frac12}\otimes Y_{3/2-n},
\end{eqnarray*}
we can replace $\varphi$ by $\varphi-u_{inn}$ and suppose that
\begin{eqnarray}
\varphi(L_1)=a_{1, 0}Y_{-\frac12}\!\otimes\! Y_{\frac32}+a_{1,
2}Y_{\frac32}\!\otimes\! Y_{-\frac12}\; (\hskip-10pt\mod {\mathcal L}\ot {\mathcal L}).
\end{eqnarray}

Applying $\varphi$ to $[L_{-1}, L_1]=2L_0$, we have $a_{-1, i}=0$ if
$i\ne 0, \pm1$, and $2a_{-1,-1}\!+\!a_{-1,0}\!+\!2a_{1,0}
=a_{-1,0}\!+\!2a_{-1,1}\!+\!2a_{1,2}=0$.

Applying $\varphi$ to $[\,L_{2},L_{-1}]=-3L_{1}$,
$[\,L_{-2},L_{1}]=3L_{-1}$ and $[L_2, L_{-2}]=-4L_0$, we obtain
$a_{\pm1,i}=0$ for all $i\in\z$, $a_{2, p}=0$ if $p\ne 0, 1, 2, 3$
$a_{-2, q}=0$ if $q\ne -2, -1, 0, 1$, and $a_{2,2}=-a_{2, 1}=3a_{2,
0}$, $a_{2, 3}=-a_{2, 0}=a_{-2, -2}=-a_{-2, 1}$, $a_{-2,
-1}=-3a_{-2, -2}=-a_{-2, 0}$.

Moreover, replaced by $\varphi$ by suitable multiples of $v_{inn}$
for $v=Y_{\frac12}\otimes Y_{-\frac12}-Y_{-\frac12}\otimes
Y_{\frac12}$ (here $v\cdot L_{\pm1}=0$),  one can suppose that
$a_{\pm2, i}=0$ for all $i\in\z$.

Since all $L_n, n\in\z$ are generated by $L_{\pm1}, L_{\pm 2}$, then
we have $a_{n, i}=0$ for all $i\in\z$. However, by $[L_{-n},
I_n]=nI_0\in\mathcal C$, we have $[L_{-n},
\varphi(I_n)]=0(\hskip-5pt\mod {\mathcal L}\ot {\mathcal L})$.
Then we have $\varphi(I_{n})=0$ if $n$ is even, and
$\varphi(I_{n})=b_n(Y_{-\frac n2}\ot
Y_{\frac{3n}2}-Y_{\frac{3n}2}\ot Y_{-\frac n2})$ if $n$ is odd.

Moreover, applying $\varphi$ to $[L_m, I_n]=nI_{m+n}$ we have
$\varphi(I_n)=0(\hskip-5pt\mod {\mathcal L}\ot {\mathcal L})$.

 Now, according to the proof of Theorem \ref{theo}, for any
$0\ne n\in\z$, we can replace $\varphi$ by $\varphi-\sg$ for some
$\sg\in{\mathcal D}_s$, we can suppose that  and suppose that
\begin{eqnarray*}
\varphi(L_n)\!&=&0,\\
\varphi(I_n)\!&=&\!0.
\end{eqnarray*}

It is easily to prove that $\varphi(Y_r)=0$.

The proofs of (ii) and (iii) are same as that in \cite{HLS}. \QED

Denote by $Y=\C\{Y_{r}\mid r\in\z+\frac12\}$, then $Y$ is a
${\mathcal L}$-module. From the proof of the Theorem \ref{theo4}, we
obtain the following result, which is very useful in determining Lie
bialgebra (super-bialgebra) structures on some Lie (super) algebras,
for example, the $N=2$ superconformal  Neveu-Schwarz algebra, etc..

\begin{coro}\label{coro45}
$H^1(\frak v, Y\ot Y)=H^1(\mathcal L, Y\ot Y)=0$. \qed
\end{coro}

\begin{rema}

We can also use the above ideas to determine Lie bialgebra
structures on some other Lie algebras and Lie superalgebras related
to the twisted Heisenberg-Virasoro algebra, for example, the twisted
Schr\"{o}dinger-Virasoro Lie algebra $($for which, the direct
calculation is very complicated, see \cite{FLL}$)$, the $N=2$
superconformal algebra $($\cite{LCZ}$)$, etc..
\end{rema}

\vskip30pt \centerline{\bf ACKNOWLEDGMENTS}

This work was partially supported by the NNSF (No. 11071068, No. 11101285), the ZJNSF(No.
D7080080,Y6100148), the "New Century 151 Talent Project" (2008),
the "Innovation Team Foundation of the Department of Education" (No.
T200924) of Zhejiang Province,  the
Shanghai Natural Science Foundation(11ZR1425900) and the Innovation
Program of Shanghai Municipal Education Commission (11YZ85).

The part of works was done by the authors visit in Kavli Institute for Theoretical Physics China(KITPC). They would like to express his special gratitude to
KITPC for financial support.

\end{document}